\documentclass[12pt]{article}
\usepackage{stmaryrd}
\usepackage{bbm}
\usepackage[french,english]{babel}

\usepackage{setspace,mathtools,amsmath,amssymb,amsthm,fullpage,outline,centernot,slashed} 
\usepackage[pdftex]{graphicx}
\usepackage{mathrsfs} 
\usepackage{enumerate}
\usepackage{accents}
\usepackage[colorlinks=true,citecolor=black,linkcolor=black,urlcolor=blue]{hyperref}
\usepackage{upgreek}

\theoremstyle{plain}
\newtheorem{theorem}{Theorem}
\newtheorem{lemma}[theorem]{Lemma}

\theoremstyle{definition}

\theoremstyle{remark}

\newcommand{\R}{\mathbb{R}}

\newcommand{\T}{\mathbb{T}}

\newcommand{\half}{\frac{1}{2}}
\newcommand{\wo}{\backslash}

\newcommand{\ben}{\begin{enumerate}}
\newcommand{\een}{\end{enumerate}}
\newcommand{\bit}{\begin{itemize}}
\newcommand{\eit}{\end{itemize}}
\DeclareMathOperator{\supp}{supp}
\def\bal#1\eal{\begin{align*}#1\end{align*}}

\DeclareMathSymbol{\mlq}{\mathord}{operators}{``}
\DeclareMathSymbol{\mrq}{\mathord}{operators}{`'}

\def\XXint#1#2#3{{\setbox0=\hbox{$#1{#2#3}{\int}$ }
\vcenter{\hbox{$#2#3$ }}\kern-.6\wd0}}

\title{From point vortices to vortex patches in self-similar expanding configurations}
\author{Samuel Zbarsky, Princeton University}

\begin{document}
\maketitle
\begin{abstract}
The main result is that given a generic self-similarly expanding configuration of 3 point vortices that start sufficiently far out, we can instead take compactly supported vorticity functions, and the resulting solution to 2D incompressible Euler will evolve like a nearby point vortex configuration for all time, with the size of the patches growing at most as $t^{1/4+\epsilon}$ and the distance between them growing as $\sqrt{t}$.
\end{abstract}
\section{Introduction}
We study the 2D Euler equation for $u:\R^2\to\R^2$
\[
\partial_tu+u\cdot\nabla u+\nabla p=0, \nabla\cdot u=0
\]
The equation can be rewritten in terms of the vorticity
\[
\omega=\partial_1u_2-\partial_2u_1
\]
as follows:
\bal
\partial_t\omega(x)&=u(x)\cdot\nabla\omega(x)\\
u&=K*\omega
\eal
where, by rescaling time to avoid factors of $2\pi$, we can take
\[
K(x)=x^\perp/|x|^2.
\]
Thus the vorticity is transported by $u$, which is generated as a singular integral of the vorticity. In this paper, we construct solutions that consist of three vortex patches (not necessarily smooth), each growing slowly in time, with the distance between the patches growing as $\sqrt{t}$. We will obtain that the trajectories of the centers of mass of these patches behave approximately like point vortices, so we will first discuss what is known about point vortex systems. A point vortex system consists of $n$ point vortices, with masses $\Omega_i$ and positions $x_i$, whose motion is described by the ODE system
\[
\frac{d}{dt}x_i=\sum_{j\ne i}\Omega_j\frac{(x_i-x_j)^\perp}{|x_i-x_j|^2}.
\]
Here, as with the vorticity formulation of 2D Euler, we rescaled time to avoid factors of $2\pi$. This ODE is meant to model a fluid in which vorticity is highly concentrated around a few points. Some information about the behavior of solutions to this ODE can be found in \cite{Aref07}. While there are specific solutions in which vortices collide, for generic initial data, this does not happen~\cite{DurrPulvurenti82,MarchioroPulvurenti94}. Rigorous justification of the point vortex model is provided in~\cite{MarchioroPulvurenti93}. They show that if one replaces point vortices by signed $L^\infty$ localized vorticity, the solution to Euler will approximate the solution to the ODE over a fixed time interval. The assumptions on the $L^\infty$ bound of the solution are then significantly weakened in independent and simultaneous works by Marchioro in~\cite{Marchioro98} and by Serfati in~\cite{Serfati98}, with~\cite{Serfati98} giving better approximation of point vortices. One can view these as being almost $L^1$ results. Both the assumptions on the $L^\infty$ bound on vorticity and the conclusion are further improved by Serfati in~\cite{OtherSerfati98}.

There are several observations about the long-term behavior of solutions to point vortex systems. First, if all masses are positive, then the solution will remain bounded. 
Second, it is easy to obtain solutions where two point vortices with masses $\Omega$ and $-\Omega$ go off to infinity with their velocity approaching some nonzero limit. 
Third, there are point vortex systems that expand and spiral in a self-similar way so that the distance between the point vortices grows as $\sqrt{t}$. An analysis of such self-similarly evolving 3-vortex system can be found in \cite{Aref07}. Some self-similarly evolving 4 and 5 vortex systems are constructed and analyzed in \cite{NovikovSedov79}. Some numerics for self-similarly evolving systems with more vortices may be found in \cite{Kudela14}. 

The centers of mass of the patches we construct will behave approximately like a self-similarly expanding 3-vortex system, see Figure 1 (for the precise statement, see Theorem~\ref{thm:threepatches} below).
\begin{figure}[h!]
\centering
\includegraphics[scale=0.45]{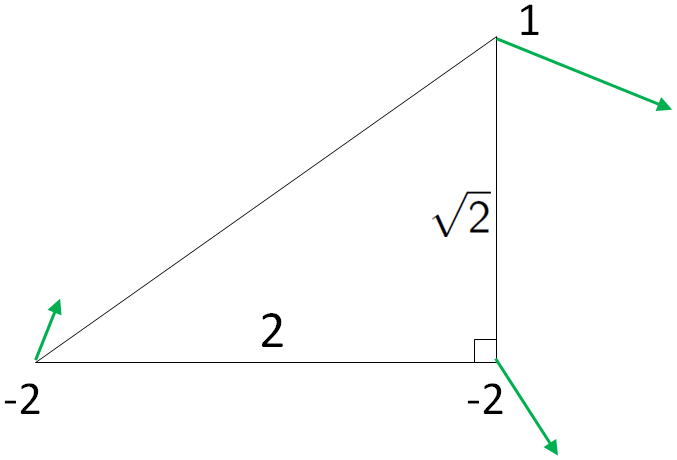}
\includegraphics[scale=0.45]{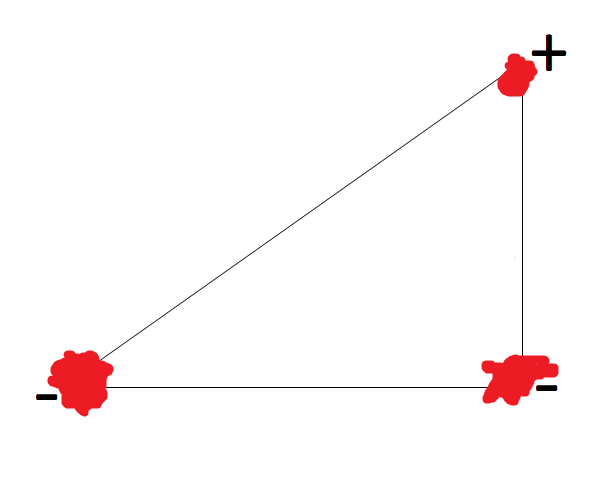}
\caption{An example of a self-similarly expanding 3-vortex system and a corresponding solution involving three vortex patches}
\end{figure}
In fact, part of the result may be viewed as an orbital stability statement for the centers of mass of the three patches relative a self-similarly expanding system. The radius of each of the three patches will grow at most like $t^{1/4+\epsilon}$. We will get this behavior for any three initial vortex patches, provided they have the right signs, are initially supported inside specified open balls, and have the correct total vorticity. In particular, this is the first construction where the support of the vorticity is known to go off to infinity that does not rely on symmetry. One consequence of not relying on symmetry is that this proof may generalize to exterior domains, where the shape of the domain itself may have no symmetries.

Before understanding the behavior of the three vortex patches, we first need to understand a single vortex patch, so we will now discuss the previously known results regarding vortex patches.
Yudovich \cite{Yudovich63} showed global well-posedness for solutions with $\omega\in L^1\cap L^\infty$. Given global well-posedness, it is natural to study the long-term behavior of vorticity. It was shown by Kirchhoff that elliptical patches will rotate uniformly~\cite{Kirchoff1874}. Other rotating solutions with $m$-fold symmetry bifurcating from the disk, called V-states, were found numerically by Deem--Zabusky~\cite{DeemZabusky78} and proved to exist by Burbea~\cite{Burbea82}. For other results about rotating solutions, see \cite{GomezSerranoetal19} and results they cite. Aside from such special solutions, it is known that if the vorticity is the indicator function of a set with $C^{k,\gamma}$ boundary, then this regularity of the boundary will continue for all time, as shown by Chemin~\cite{Chemin91}, Bertozzi--Constantin~\cite{ConstantinBertozzi93}, and Serfati~\cite{Serfati94}. For other results concerning regularity and long-term behavior of vortex patches, see \cite{ElgindiJeong19} and results they cite.



However, very little is known if no additional regularity is assumed. In particular, one can ask what happens if vorticity is initially compactly supported and $L^\infty$. We will go over some results bounding the expansion of the support, known as vorticity confinement results. It is easy to see that the radius of the support can grow at most linearly, since $u$ is bounded. If vorticity is of a definite sign, then the radius of the support grows more slowly. In fact, it is not known whether it ever goes to infinity. Marchioro~\cite{Marchioro94} showed an upper bound of $t^{1/3}$ by using conservation of the second moment of vorticity. This was independently improved to $(t\log t)^{1/4}$ by Iftimie--Sideris--Gamblin~\cite{InftimieSiderisGamblin99} and to $t^{1/4}\log \circ\cdots\circ\log t$ by Serfati \cite{Serfati98}, with the improvement coming largely from using conservation of the center of mass of the vorticity. Compare this with the present work, in which we get the slightly worse bound of $t^{1/4+\epsilon}$ for each of the vortex patches. There are also several other vorticity confinement results, including getting similar confinement bounds, but depending on the $L^{2+\epsilon}$ norm rather than the $L^\infty$ norm of vorticity~\cite{LopesLopes98}. Compare this to the result in \cite{Serfati98}, which requires an $L^\infty$ bound, but the constant in the confinement bound has very weak dependence on the  $L^\infty$ norm (it is linear in $\log \circ\cdots\circ\log ||\omega||_{L^\infty})$, depending mostly on the $L^1$ norm. There are also various bounds on confinement of positive compactly supported vorticity in other domains. In particular, on the upper half-plane, the $x$ coordinate of the center of mass of vorticity is at least $ct$ and the $y$ coordinate of points in the support is bounded by $C(t\log t)^{1/3}$~\cite{Iftimie99}, while the $x$ coordinate of  points in the support is at least $-C(t\log t)^{1/2}$ \cite{SecondIftimieetal03}. The latter work also analyzes what possible weak limits a positive vorticity solution can have on a half-plane (under appropriate rescaling). In exterior domains, the radius of the support is bounded by $Ct^{1/2}$, with further improvements in the exponent when the domain is the exterior of a disk~\cite{Marchioro96,Iftimieetal07}. On $\T\times\R$, the $y$ coordinate of points in the support is bounded by $Ct^{1/3}\log^2 t$~\cite{ChoiDenisov19}. A survey of various related results can be found in \cite{Iftimie07}. However, things are very different if you allow mixed-sign vorticity. \cite{InftimieSiderisGamblin99} contains a construction of a compactly supported  positive vorticity vortex patch in the first quadrant of the plane, reflected into the other quadrants with changing sign, whose support grows linearly. However, the proof relies heavily on the symmetry, so it is very unstable and can only give this result for a system with total vorticity 0.

Returning to the plane, but now without a definite sign, there is a vorticity confinement result by Iftimie--Lopes Filho--Nussenzveig Lopes that addresses the question of weak limits under appropriate rescaling~\cite{Iftimieetal03}. This result states that if we define
\[
\tilde\omega_\alpha(x,t)=t^{2\alpha}\omega(t^\alpha x,t)
\]
for $\alpha>1/2$, then
\[
\tilde\omega_\alpha\xrightharpoonup[t\to\infty]{} m\delta_0
\]
in the weak-* sense for measures where $m=\int\omega dx$ and $\delta_0$ is the Dirac delta. The authors interpret this as showing confinement of net vorticity to a radius of $\sqrt{t}$, but this result still allows for strange examples like having both positive and negative vortex patches moving away from the origin in different directions and being at distance $t^{2/3}\log t$.  Our result shows that we cannot take $\alpha=1/2$  in the statement of~\cite{Iftimieetal03} and thus, in a sense, net vorticity is moving off to infinity. In fact, the solution given here should, modulo a rotation by a logarithmically growing angle, weakly converge to a sum of three delta masses under the rescaling with $\alpha=1/2$.

The last previous result we discuss is a paper by Iftimie--Marchioro~\cite{IftimieMarchioro18}, which looked at a toy model of the construction given here and showed confinement of the vortex patches. The toy model consists of taking a self-similarly expanding point vortex system, replacing only one of the point vortices with a patch, assuming that the trajectory of the other point vortices is fixed, and seeing how the patch evolves. The purpose of looking at the toy model was to sidestep the issue of stability for self-similar point vortex systems. and only worry about confinement of vorticity. For some configurations, they bound the radius of the support of the patch as growing no faster than $t^{(1+\alpha)/3}$ for $\alpha$ some constant that depends on the configuration, is always positive, and is less than 1/2. This means that the vortex patch size grows slower than the distance between patches. In this paper, we solve the full problem for self-similar expanding three-vortex systems, replacing each of the three vortices with vortex patches, allowing a generic self-similarly expanding configuration of 3 vortices, and obtaining that the radius of support grows at most as $t^{1/4+\epsilon}$ for any given $\epsilon>0$. The improvement in the exponent from $1/3$ to $1/4$ comes from actually analyzing the stability and keeping track of the center of mass, in the same way that using conservation of center of mass gives the same improvement for a single vortex patch. To get rid of $\alpha$, we obtain a better bound on moment growth by noticing that  one of the expressions that shows up in the expression for moment derivatives is the approximate derivative of another expression and thus obtaining some cancellation in the most troublesome term (one can also think of this as a renormalization of the moments). Our result is limited to 3 vortex system due to the stability analysis, and if one found stably growing systems of 4 or more vortices, the confinement result would most likely carry over with little modification. However, one needs sufficiently good stability results; orbital linear stability, which may not be hard to obtain for some systems, is not enough. Our assumptions on the patches, same as in \cite{IftimieMarchioro18}, will be that the vorticity $\omega$ is compactly supported and  $L^{2+\epsilon}$ (note that a smaller $\epsilon$ both expands this class and improves vortex patch confinement, so this is not a tradeoff). Technically, we need to assume that $\omega\in L^\infty$ to make use of the well-posedness theory, but all constants in the proof will only depend on the  $L^{2+\epsilon}$ bound. We need this bound in many places in the proofs in order to bound integrals of $\frac{1}{|x-a|}\omega(x)$ using Holder's inequality.

\section{Point vortex systems}
In order to state the main result, we first need to analyze some properties of expanding systems of three point vortices. Take any three vortex system. It has four conserved quantities, each of which is easy to check.
\begin{enumerate}
\item $X=\sum_i x_i\Omega_i$ (this has two components)
\item $I=\sum_i \Omega_i|x_i|^2$
\item $E=\sum_{i\ne j} \Omega_i\Omega_j\log|x_i-x_j|$.
\end{enumerate}
Now suppose we take a self-similarly expanding solution with three point vortices and nonzero total mass. By moving the origin, we can assume that $X=0$. Then self-similarity ensures that $I=Ct$ for some constant $C$.Then conservation of $I$ and $E$ ensure that 
\begin{align}
\Omega_1\Omega_2+\Omega_1\Omega_3+\Omega_2\Omega_3&=0\label{eq:cond1}\\
I&=0\label{eq:cond2}
\end{align}
We move the rest of the analysis to the following lemma, proved in the appendix:
\begin{lemma}\label{configlemma}
Suppose we have a point vortex system $(y_i,\Omega_i)$ satisfying $X=0$, \eqref{eq:cond1}, and \eqref{eq:cond2}. Suppose furthermore that $y_1,y_2,y_3$ are not the vertices of an equilateral triangle and are not collinear and that they evolve as a self-similarly expanding system.
\begin{enumerate}
\item $\sum\Omega_i\ne 0$.
\item Take the subspace $V$ of vortex locations that satisfy $X=0$. There exists a  two-dimensional surface $S$ through $(y_i)$ and coordinates on some neighborhood $U$ of $(y_1,y_2,y_3)$ in $A$ such that two of the coordinates are the angle of rotation and the ratio of dilation needed to hit some $z\in S$ and the other two coordinates are $I(z)$ and $E(z)$ at the resulting point of $S$.
\end{enumerate}
\end{lemma}
Note that for the second part of the lemma statement, it is important that $I$ is evaluated at the point $z\in S$, not at the original point $x$.

The conditions in the lemma statement are generic for self-similarly expanding 3-vortex configurations; one system satisfying the hypotheses of the lemma is the following example, taken (up to sign reversal) from \cite{MarchioroPulvurenti94}. Let $\Omega_1=-2$, $\Omega_2=-2$, $\Omega_3=1$ and $x_1=(-1,0)$, $x_2=(1,0)$, $x_3=(1,\sqrt{2})$. Then translate to achieve $X=0$. This example is shown in Figure 1.
\section{Result statement and stability of centers of mass}
We can now state the main result precisely.
\begin{theorem}\label{thm:threepatches}
Take an arrangement of three points $(y_i,\Omega_i)$ satisfying the conditions of Lemma~\ref{configlemma}. Let $\epsilon>0$, be arbitrary and small, $M>0$, $\rho>0$ be arbitrary. Then there exists some $T>0$ so that for $t_0>T$, if we take the solution $\{(\sqrt{t_0}y_i,\Omega_i)\}$ and replace each point vortex $(x_i,\Omega_i)$ with an $L^\infty$ vorticity function $\omega_i$ such that:
\begin{enumerate}
           \item The center of mass of the whole system is still 0.\label{cond1}
	\item $\supp\omega_i\subseteq B(x_i,\rho)$.
	\item $||\omega_i||_{L^{2+\epsilon}}\le M$.
	\item $\omega_i$ has definite sign.
	\item $\int\omega_i=\Omega_i$.
\end{enumerate}
then at each later time $t$, there exists some $\{(z_i(t))\}\in S$ with $|z_i-y_i|<\epsilon$, some angle $\beta(t)$, and some factor $\gamma(t)=(1+O(\epsilon))\sqrt{t})$ such that, letting $x_i=R_{\beta}\gamma z_i$, the solution at time $t$ is $\sum \omega_i$ with:
\begin{enumerate}
          \item $x_i=\frac{1}{\Omega_i}\int \omega_i(x)x dx$.
	\item $\supp\omega_i\subseteq B(x_i,\epsilon t^{1/4+\epsilon})$.
	\item $\int\omega_i=\Omega_i$.
	\item $\omega_i$ has definite sign.
\end{enumerate}
\end{theorem}
There are a couple of comments about the statement. First, because of how we defined $x_i$ and $S$, we have that $\beta$ and $\gamma$ are uniquely defined by the arrangement. Second, condition~\ref{cond1} in the theorem statement is simply for convenience--since $\sum\Omega_i\ne 0$, we could restate the theorem without this condition, but adding a translation to move the center of mass to 0. Third, for the theorem as stated, for different $\epsilon$ one needs different data. However, by being much more careful, and with some modification of the bootstrap assumptions, one can probably modify the proof to find some $T$ so that the radius of the support of each $\omega_i$ grows like $t^{1/4+o(1)}$. Finally, a version of this theorem and proof probably holds for exterior domains. One would need to use the fact that in exterior domains, when $|x|,|y|\sim\sqrt{t}$, we have $K(x,y)\approx \frac{(x-y)^\perp}{|x-y|^2}$ (using the more precise asymptotics found in \cite{Iftimieetal07}). The theorem statement would need to be modified slightly since center of mass is no longer conserved.

\textbf{Proof}.

$k$ will denote an even integer and $\delta$ will denote some sufficiently small constant that can depend on $k$ and $R$ will denote some large constant that depends on $\delta$. At the end of the proof, we will choose $k$, $\delta$, and $R$ depending on the initial configuration of point vortices, as well as $\epsilon, \rho, M$. We will then choose $T$ large enough depending on $\epsilon, \rho, M,k,\delta$. All constants $C$ in the statement and proof (including implicit constants hidden by $O$ notation) can depend on the initial configuration of point vortices, as well as $\epsilon, \rho, M,k$, but not on $R$ or $\delta$. The letters $C$ may be used for different constants on different lines. We will use $\hat O$ if we're allowing the implicit constant to depend on the initial configuration of point vortices and on nothing else.
At any time $t$, let
\bal
x_i(t)&=\frac{1}{\Omega}\int x\omega_i(x)dx\\
I_{k,i}(t)&=\int |x-\tilde x_i|^k\omega_i(x)dx
\eal
We will have the following bootstrap assumptions:
\begin{enumerate}
\item $x_i=R_{\beta}\gamma z_i$ with $|z_i-y_i|<\epsilon^2$ for some angle $\beta(t)$, and some factor $\gamma(t)$  satisfying $|\frac{\gamma}{\sqrt{t}}-1|<\epsilon$ \label{bootstrapxi}
\item $I_{2,i}< t^{\epsilon/2}$\label{bootstrapI2i}
\item $I_{k,i}< t^{k(1+\epsilon)/4}$\label{bootstrapIki}
\item $\omega_1, \omega_2, \omega_3$ are three $L^{2+\epsilon}$ compactly supported functions of definite sign, $\|\omega_i\|_{L^{2+\epsilon}}\le M$
\item $\supp\omega_i\subseteq B(x_i,\epsilon t^{1/4+\epsilon})$, that is $|p-x_i|<\epsilon t^{1/4+\epsilon}$ for any $p\in\supp\omega_i$.\label{bootstrapsupport}
\end{enumerate}
These assumptions hold at time $t_0$ as long as $T$ is big enough. If they always hold, we are done, so we can assume that the first time when one of them fails is $T_*<\infty$.

First, we want to understand the ODE satisfied by the triple $(x_1,x_2,x_3)$ of centers of mass of the patches in order to verify bootstrap assumption~\ref{bootstrapxi}. First note that from the conservation of the center of mass of the vorticity, we get that $\sum \Omega_ix_i=0$, so the center of mass of $z_i$ stays at 0. We will use the notation $I_x=\sum \Omega_i|x_i|^2$ and $I_z=\sum \Omega_i|z_i|^2$. Because of condition~\ref{cond1}, there is so need to make this distinction for $E$. For this calculation, we note that from the bootstrap assumptions, for $x\in\supp \omega_i$, $y\in\supp \omega_j$ with $j\ne i$, we have
\[
\frac{(x-y)^\perp}{|x-y|^2}=\frac{(x_i-x_j)^\perp}{|x_i-x_j|^2}+A_1(x-x_i)+A_2(y-x_j)+O\left(\frac{|x-x_i|^2+|y-x_j|^2}{|x_i-x_j|^3}\right)
\]
where $A_1,A_2$ are some linear functions dependent on $t$. Then
\bal
\frac{d}{dt}x_i&=\frac{1}{\Omega_i}\iint \omega_i(x)\omega_i(y)\frac{(x-y)^\perp}{|x-y|^2}  dxdy+\sum_{j\ne i} \iint \omega_i(x)\omega_j(y)\frac{(x-y)^\perp}{|x-y|^2}  dxdy\\
&=\frac{1}{\Omega_i}\Bigg[\half\iint \omega_i(x)\omega_i(y)\left(\frac{(x-y)^\perp}{|x-y|^2}+\frac{(y-x)^\perp}{|y-x|^2}\right)  dxdy\\
&\qquad+\sum_{j\ne i}\Omega_i(x)\Omega_j(y)\frac{(x_i-x_j)^\perp}{|x_i-x_j|^2}+O\left(\sum_{j=1}^3\frac{ I_j}{t^{3/2}}\right)\Bigg]\\
&=\sum_{j\ne i}\Omega_j(y)\frac{(x_i-x_j)^\perp}{|x_i-x_j|^2}+O(t^{\epsilon/2-3/2}).
\eal
This means that the system $(x_i)$ evolves as point vortices, up to some error. We now look at the nearly-conserved quantities:
\bal
\left|\frac{dI_x}{dt}\right|&=\left|\frac{d}{dt}\sum_i \Omega_i|x_i|^2\right|=\left|\sum_i 2\Omega_ix_i\cdot\left(\sum_{j\ne i}\Omega_j(y)\frac{(x_i-x_j)^\perp}{|x_i-x_j|^2}+O(t^{\epsilon/2-3/2})\right)\right|\\
&=O(t^{\epsilon/2-1})
\eal
where the fist terms canceled algebraically (this being the same calculation that gave us the conserved quantity in the first place) and where we used bootstrap assumption~\ref{bootstrapxi} to get that $|x_i|=O(\sqrt{t})$. Furthermore, at time $t_0$, we have $I_x=O(t_0^{1/2})$. Integrating in $t$, we get $I_x=O(t_0^{1/2}+t^{\epsilon/2})$.
This then gives us that
\[
I_z=O(I_x/t)=O\left(\frac{t_0^{1/2}}{t}+t^{\epsilon/2-1}\right)=O(t^{-1/2}).
\]
Finally,
\bal
\left|\frac{dE}{dt}\right|&=\left|\frac{d}{dt}\sum_{i\ne j} \Omega_i\Omega_j\log|x_i-x_j|\right|=\sum_{i\ne j}\frac{ \Omega_i\Omega_j}{|x_i-x_j|} O(t^{\epsilon/2-3/2})
\\
&=O(t^{\epsilon/2-2})
\eal
where once again the non-error terms cancel precisely. Integrating in $t$, we get that $|E(t)-E(t_0)|= O(t_0^{\epsilon/2-1})$. Therefore, by choosing $T$ sufficiently large, we can guarantee that both $I_z$ and $E$ change by small amounts from their values at $y_i$, so $|z_i-y_i|<\epsilon^2$, so that part of bootstrap assumption~\ref{bootstrapxi} is maintained.

Now we note that
\begin{align*}
\frac{d}{dt}|x_1|^2&=2x_1\cdot \sum_{j\ne 1}\Omega_j(y)\frac{(x_1-x_j)^\perp}{|x_1-x_j|^2}+O(t^{\epsilon/2-1})=2z_1\cdot\sum_{j\ne 1}\Omega_j(y)\frac{(z_1-z_j)^\perp}{|z_1-z_j|^2}+O(t^{\epsilon/2-1})\\
&=2y_1\cdot\sum_{j\ne 1}\Omega_j(y)\frac{(y_1-y_j)^\perp}{|y_1-y_j|^2}+\hat O(\epsilon^2)+O(t^{\epsilon/2-1})=|y_1|^2+\hat O(\epsilon^2)+O(t^{\epsilon/2-1})
\end{align*}
where the last equality is because of the expansion rate of the point vortex system. From this, we get that
\[
\left|\gamma(t)^2-t\right|=\left|\frac{|x_1|^2}{|z_1|^2}-t\right|=\left|\frac{|x_1|^2}{|y_1|^2}-t-\hat O(t\epsilon^2)\right|=\hat O(\epsilon^2t)+O(t^{\epsilon/2-1}).
\]
From this, and assuming that $\epsilon$ is sufficiently small while $T$ is sufficiently large, we get that
\[
\left|\frac{\gamma}{\sqrt{t}}-1\right|<\epsilon
\]
which is the last remaining part of bootstrap assumption~\ref{bootstrapxi}. Bootstrap assumption 4 is simply a consequence of the vorticity being transported by an incompressible flow (generated by divergence-free vector field $u$).
For the other bootstrap assumptions, there are two cases: $T_*<t_0+t_0^{9/10}$ and $T_*\ge t_0+t_0^{9/10}$. We will handle the latter case first, since the former case uses weaker versions of estimates that we'll need to derive along the way.

\section{Long time behavior}\label{sec:longterm}
In this section, we assume that $T_*\ge t_0+t_0^{9/10}$.

We first prove that  bootstrap assumptions~\ref{bootstrapI2i} and \ref{bootstrapIki} are maintained by bounding $\frac{d}{dt}I_{k,i}$. We will mostly treat them together, as many of the calculations can be done for any $k$, and we will simply plug in 2 for $k$ when we need to.
For definiteness, we will take $i=1$ below, and we assume that $\omega_1\ge 0$.  The proof for $\omega_2$ and $\omega_3$ and for negative vorticity is identical.
We want to bound the growth of $I_{k,1}$.  Let $v_j$ be the velocity field generated by $\omega_j$. Then for $x\in \supp(\omega_1)$, we Taylor expand $K$ to get that for some linear function $A_t$,
\begin{align}
v_2(x)&=\frac{(x_1-x_2)^\perp}{|x_1-x_2|^2}+\int A_t(y-x_2)  \omega_2(y)dy\nonumber\\
&\qquad-(x-x_1)\cdot(x_1-x_2)\frac{(x_1-x_2)^\perp}{|x_1-x_2|^4}-(x-x_1)\cdot(x_1-x_2)^\perp\frac{x_1-x_2}{|x_1-x_2|^4}\nonumber\\
&\qquad+O\left(\frac{|x-x_1|^2}{|x_1-x_2|^3}\right)+O\left(\int\frac{|y-x_1|^2}{|x_1-x_2|^3}\omega_2(y)dy\right)\nonumber\\
&=\frac{(x_1-x_2)^\perp}{|x_1-x_2|^2}-(x-x_1)\cdot(x_1-x_2)\frac{(x_1-x_2)^\perp}{|x_1-x_2|^4}-(x-x_1)\cdot(x_1-x_2)^\perp\frac{x_1-x_2}{|x_1-x_2|^4}\nonumber\\
&\qquad+O\left(\frac{|x-x_1|^2}{|x_1-x_2|^3}\right)+O\left(\frac{I_{2,2}}{|x_1-x_2|^3}\right)\label{eq:v2}
\end{align}
where we used that $x_2$ is the center of mass of $\omega_2$ to eliminate one of the terms. We have a similar expression for $v_3(x)$. Then
\begin{equation}\label{eq:dx1dt}
\frac{dx_1}{dt}=\int (v_2(x)+v_3(x))\omega_1(x)dx=v_2(x_1)+v_3(x_1)+\tilde v
\end{equation}
where
\[
\tilde v= O\left(\frac{I_{2,1}+I_{2,2}+I_{2,3}}{t^{3/2}}\right)=O(t^{\epsilon/2-3/2}).
\]
We then have
\begin{align}
\frac{d}{dt}I_{k,1}&=\frac{d}{dt}\int |x-x_1|^k\omega_1(x)dx\nonumber\\
&=\iint k|x-x_1|^{k-2}\frac{(x-x_1)\cdot(x-y)^{\perp}}{|x-y|^2}\omega_1(x)\omega_1(y)dxdy\nonumber\\
&\qquad+\int k|x-x_1|^{k-2}(x-x_1)\cdot (v_2(x)-v_2(x_1)+v_3(x)-v_3(x_1))\omega_1(x)dx\label{eq:ddtIk1}.
\end{align}

We first deal with the first term in the same way that it is done in $\cite{InftimieSiderisGamblin99}$. First, we note that in the special case $k=2$, we symmetrize in $x$ and $y$, and that term vanishes. When $k>2$ even, we use the fact that $x_1$ is the center of mass to subtract 0. Note that because $|(x-x_1)\cdot(x_1-y)^{\perp}|=|(x-y)\cdot(x_1-y)^{\perp}|\le |x-y||x_1-y|$, all terms in the expressions below are absolutely integrable, so the rearrangements and splitting are all valid.
\begin{align}
\iint &|x-x_1|^{k-2}\frac{(x-x_1)\cdot(x-y)^{\perp}}{|x-y|^2}\omega_1(x)\omega_1(y)dxdy=\nonumber\\
&=\iint |x-x_1|^{k-2}\frac{(x-x_1)\cdot(x_1-y)^{\perp}}{|x-y|^2}-|x-x_1|^{k-2}\frac{(x-x_1)\cdot(x_1-y)^{\perp}}{|x-x_1|^2}\omega_1(x)\omega_1(y)dxdy\nonumber\\
&=\iint_{S_1} |x-x_1|^{k-2}\frac{(x-x_1)\cdot(x_1-y)^{\perp}}{|x-y|^2}-|x-x_1|^{k-2}\frac{(x-x_1)\cdot(x_1-y)^{\perp}}{|x-x_1|^2}\omega_1(x)\omega_1(y)dxdy\nonumber\\
&\qquad+\iint_{S_2}|x-x_1|^{k-2}\frac{(x-x_1)\cdot(x_1-y)^{\perp}}{|x-y|^2}-|x-x_1|^{k-2}\frac{(x-x_1)\cdot(x_1-y)^{\perp}}{|x-x_1|^2}\omega_1(x)\omega_1(y)dxdy\nonumber\\
&\qquad+\iint_{S_3}|x-x_1|^{k-2}\frac{(x-x_1)\cdot(x_1-y)^{\perp}}{|x-y|^2}-|x-x_1|^{k-2}\frac{(x-x_1)\cdot(x_1-y)^{\perp}}{|x-x_1|^2}\omega_1(x)\omega_1(y)dxdy\label{eq:threeregions}
\end{align}
where
\begin{align*}
S_1&=\{|x-x_1|<|y-x_1|/2\}\\
S_2&=\{|y-x_1|/2\le |x-x_1|\le 2|y-x_1|\}\\
S_3&=\{|x-x_1|>2|y-x_1|\}
\end{align*}
For the first term in $\eqref{eq:threeregions}$, we use
\begin{align}
&\left|\iint_{S_1} |x-x_1|^{k-2}\frac{(x-x_1)\cdot(x_1-y)^{\perp}}{|x-y|^2}-|x-x_1|^{k-2}\frac{(x-x_1)\cdot(x_1-y)^{\perp}}{|x-x_1|^2}\omega_1(x)\omega_1(y)dxdy\right|\nonumber\\
&\qquad\le\iint_{S_1} 4|x-x_1|^{k-2}\omega_1(x)\omega_1(y)dxdy\nonumber\\
&\qquad\le \iint_{S_1} |x-x_1|^{k-4}|y-x_1|^2\omega_1(x)\omega_1(y)dxdy\nonumber\\
&\qquad\le CI_{k,1}^{\frac{k-4}{k}}I_{2,1}\le Ct^{\epsilon/2}I_{k,1}^{\frac{k-4}{k}}\label{S1}
\end{align}
For the second term in $\eqref{eq:threeregions}$, we symmetrize in $x$ and $y$ and get
\begin{align}
&\left|\iint_{S_2}|x-x_1|^{k-2}\frac{(x-x_1)\cdot(x_1-y)^{\perp}}{|x-y|^2}-|x-x_1|^{k-2}\frac{(x-x_1)\cdot(x_1-y)^{\perp}}{|x-x_1|^2}\omega_1(x)\omega_1(y)\right|\nonumber\\
&\qquad\le\left|\iint_{S_2}(|x-x_1|^{k-2}-|y-x_1|^{k-2})\frac{(x-x_1)\cdot(x_1-y)^{\perp}}{|x-y|^2}\omega_1(x)\omega_1(y)dxdy\right|\nonumber\\
&\qquad\qquad+\left|\iint_{S_2}|x-x_1|^{k-4}(x-x_1)\cdot(x_1-y)^{\perp}\omega_1(x)\omega_1(y)dxdy\right|\nonumber\\
&\qquad\le\Bigg|\iint_{S_2}(k/4-1/2)\big(|x-x_1|^{k-4}+|y-x_1|^{k-4}\big)\big(|x-x_1|^2-|y-y_1|^2\big)\nonumber\\
&\qquad\qquad\times\frac{(x-x_1)\cdot(x-y)^{\perp}}{|x-y|^2}\omega_1(x)\omega_1(y)dxdy\Bigg|\nonumber\\
&\qquad\qquad+2\left|\iint_{S_2}|x-x_1|^{k-4}|y-x_1|^2\omega_1(x)\omega_1(y)dxdy\right|\nonumber\\
&\qquad\le C\left|\iint_{S_2}|x-x_1|^{k-4}|y-k_1|^2\omega_1(x)\omega_1(y)dxdy\right|\nonumber\\
&\qquad\qquad+2\left|\iint_{S_2}|x-x_1|^{k-4}|y-x_1|^2\omega_1(x)\omega_1(y)dxdy\right|\nonumber\\
&\qquad\le CI_{k,1}^{\frac{k-4}{k}}I_{2,1}\le Ct^{\epsilon/2}I_{k,1}^{\frac{k-4}{k}}.\label{S2}
\end{align}
For the third term in $\eqref{eq:threeregions}$, we take advantage of the term we subtracted off to get
\begin{align}
&\left|\iint_{S_3} |x-x_1|^{k-2}\frac{(x-x_1)\cdot(x_1-y)^{\perp}}{|x-y|^2}-|x-x_1|^{k-2}\frac{(x-x_1)\cdot(x_1-y)^{\perp}}{|x-x_1|^2}\omega_1(x)\omega_1(y)dxdy\right|\nonumber\\
&\qquad\le C\iint_{S_3} |x-x_1|^{k-4}|y-x_1|^2\omega_1(x)\omega_1(y)dxdy\nonumber\\
&\qquad\le CI_{k,1}I_{2,1}\le Ct^{\epsilon/2}I_{k,1}^{\frac{k-4}{k}}.\label{S3}
\end{align}
Plugging \eqref{S1}, \eqref{S2}, \eqref{S3} into \eqref{eq:threeregions} when $k>2$ and remembering that the term goes away when $k=2$, we get that for $k\ge 2$ even
\begin{equation}\label{usualmomentderiv}
\left|\iint |x-x_1|^{k-2}\frac{(x-x_1)\cdot(x-y)^{\perp}}{|x-y|^2}\omega_1(x)\omega_1(y)dxdy\right|\le C(k-2)t^{\epsilon/2}I_{k,1}^{\frac{k-4}{k}}.
\end{equation}

To deal with the second term in \eqref{eq:ddtIk1}, we plug \eqref{eq:v2} into it to get
\begin{align*}
\int &|x-x_1|^{k-2}(x-x_1)\cdot (v_2(x)-v_2(x_1))\omega_1(x)dx=\nonumber\\
&=\int -2|x-x_1|^{k-2}(x-x_1)\cdot(x_1-x_2)^\perp\frac{(x-x_1)\cdot(x_1-x_2)}{|x_1-x_2|^4}\omega_1(x)\\
&\qquad+O\left(\frac{|x-x_1|^{k+1}}{|x_1-x_2|^3}\right)\omega_1(x)+O\left(\frac{I_{2,2}|x-x_1|^{k-1}}{|x_1-x_2|^3}\right)\omega_1(x)dx\\
&=\int -\frac{\sin(2\theta(x))|x-x_1|^{k}}{|x_1-x_2|^2}\omega_1(x)+O\left(|x-x_1|^{k}\frac{|x-x_1|}{|x_1-x_2|^3}\right)\omega_1(x)\\
&\qquad+O\left(\frac{I_{2,2}}{|x_1-x_2|^3}\right)|x-x_1|^{k-1} \omega_1(x)dx
\end{align*}
where $\theta(x)$ is the angle between $x-x_1$ and $x_1-x_2$. Now, using the bootstrap assumptions on $\sup |x-x_1|$, $I_{2,2}$, and $|x_1-x_2|$, as well as using Holder's inequality on the last term, we get
\begin{align}
\int& |x-x_1|^{k-2}(x-x_1)\cdot v_2(x)\omega_1(x)dx=\nonumber\\
&=\int -\frac{\sin(2\theta(x))|x-x_1|^{k}}{|x_1-x_2|^2}\omega_1(x)dx+O\left(t^{\epsilon-5/4}\right)I_{k,1}+O\left(t^{\frac{\epsilon-3}{2}}\right)I_{k,1}^{\frac{k-1}{k}}.\label{eq:dotv2}
\end{align}
If we were to use crude bounds for the first term of $\eqref{eq:dotv2}$, bounding the numerator by $I_{k,1}$, we would achieve vorticity confinement that is worse by some factor of $t^\alpha$, with $\alpha$ depending on the configuration of point vortices $(\Omega_i,y_i)$. In fact, for some configurations, our confinement result would be worse than $t^{1/2}$, so it wouldn't be enough to prevent the patches from interacting strongly, causing the whole proof to break down. For this reason, we want a better bound on this term. There is little hope of getting one for each time, but we note that we want to bound the expression in \eqref{eq:dotv2} because it appears in the derivative of $I_{k,1}$, so it is enough to get a better bound on its time average (as long as the time interval we are averaging over isn't too long). More precisely, let $f_{k,2}(x)=-\cos(2\theta(x))|x-x_1|^{k+2}$ and use the following estimate, which will be proved in section~\ref{derivest}:
\begin{equation}\label{eq:subsecresult}
\frac{d}{dt}\int f_{k,2}(x)\omega_1(x)dx=\Omega_1\int 2\sin(2\theta)|x-x_1|^k\omega_1(x) dx+O(\delta I_{k,1}+R^k\delta^{-k}).
\end{equation}
One way of thinking about this estimate is as a renormalization of the moments $I_{k,i}$, where we define new quantities of the form
\[
\hat I_{k,i}=I_{k,i}+C\frac{\int f_{k,2}(x)\omega_1(x)dx}{t}
\]
and bound their time derivatives. This is equivalent to the argument below, though the notation and organization are different.

To use \eqref{eq:subsecresult}, we plug it into into \eqref{eq:dotv2} to get
\begin{align}
\int |x-x_1|^{k-2}&(x-x_1)\cdot v_2(x)\omega_1(x)dx=\frac{-1}{2\Omega_1|x_1-x_2|^2}\frac{d}{dt}\int f_{k,2}(x)\omega_1(x)dx\nonumber\\
&+O\left(t^{\epsilon-5/4}I_{k,1}\right)+O\left(t^{\frac{\epsilon-3}{2}}I_{k,1}^{\frac{k-1}{k}}\right)+O\left(\frac{\delta I_{k,1}+R^k\delta^{-k}}{t}\right).\label{eq:dotv2second}
\end{align}
If we introduce $f_{k,3}$ as being entirely analogous to $f_{k,2}$, but with $x_3$ replacing $x_2$, we  can plug \eqref{usualmomentderiv} and \eqref{eq:dotv2second} and into \eqref{eq:ddtIk1} to get
\begin{align}
\frac{d}{dt}I_{k,1}&=O\left((k-2)t^{\epsilon/2}I_{k,1}^{\frac{k-4}{k}}\right)-\frac{k}{2\Omega_1|x_1-x_2|^2}\frac{d}{dt}\int f_{k,2}(x)\omega_1(x)dx\nonumber\\
&\qquad-\frac{k}{2\Omega_1|x_1-x_3|^2}\frac{d}{dt}\int f_{k,3}(x)\omega_1(x)dx+O\left(t^{\epsilon-5/4}\right)I_{k,1}\nonumber\\
&\qquad+O\left(t^{\frac{\epsilon-3}{2}}\right)I_{k,1}^{\frac{k-1}{k}}+O\left(\frac{\delta I_{k,1}+R^k\delta^{-k}}{t}\right)\label{eq:ddtIk1compiled}
\end{align}
We are now ready to confirm the bootstrap assumptions on $I_{2,1}$ and $I_{k,1}$. First, we note that
\begin{equation}\label{eq:f2bound}
\left|\int f_{k,2}(x)\omega_1(x)dx\right|\le \left(\sup_{x\in \supp\omega_1}|x-x_1|\right)^2\int|x-x_1|^{k}\omega_1(x)dx\le \epsilon^2t^{1/2+2\epsilon}I_{k,1}.
\end{equation}
We now set $k=2$ in order to prove that $I_{2,1}(T_*)< T_*^{\epsilon/2}$. The term that has a factor of $k-2$ in \eqref{eq:ddtIk1compiled} goes away, and we get
\begin{align}
\frac{d}{dt}I_{2,1}&=-\frac{2}{2\Omega_1|x_1-x_2|^2}\left(\frac{d}{dt}\int f_{2,2}(x)\omega_1(x)dx\right)-\frac{2}{2\Omega_1|x_1-x_3|^2}\left(\frac{d}{dt}\int f_{2,3}(x)\omega_1(x)dx\right)\nonumber\\
&\qquad+O\left(t^{\epsilon-5/4}I_{2,1}\right)+O\left(t^{\frac{\epsilon-3}{2}}I_{2,1}^{\frac{1}{2}}\right)+O\left(\frac{\delta I_{2,1}+R^k\delta^{-k}}{t}\right)\nonumber
\end{align}
We integrate in $t$ from $t_1=T_*-T_*^{2/3}\ge t_0$ to $T_*$ and plug in $I_{2,1}(t)\le t^{\epsilon/2}$ to get
\begin{align}
I_{2,1}(T_*)&\le  -\int_{t_1}^{T_*}\frac{2}{2\Omega_1|x_1-x_2|^2}\left(\frac{d}{dt}\int f_{2,2}(x)\omega_1(x)dx\right)+\frac{2}{2\Omega_1|x_1-x_3|^2}\left(\frac{d}{dt}\int f_{2,3}(x)\omega_1(x)dxdt\right)\nonumber\\
&\qquad+C\int_{t_1}^{T_*}t^{3\epsilon/2-5/4}+t^{3\epsilon/4-3/2}+\delta t^{\epsilon/2-1}+R^k\delta^{-k}t^{-1}dt+I_{2,1}(t_1)\nonumber
\end{align}
We now integrate by parts in $t$ and bound boundary terms using the bootstrap assumptions on $|x_1-x_2|$ as well as \eqref{eq:f2bound} (as well as the same bound for $f_{2,3}$). We also use the fact that $T\gg R,1/\delta$ to get
\begin{align}
I_{2,1}(T_*)&\le CT_*^{-1/2+2\epsilon}(I_{2,1}(t_1)+I_{2,1}(T_*))+\int_{t_1}^{T_*}\left(\frac{d}{dt}\frac{2}{2\Omega_1|x_1-x_2|^2}\right)\int f_{2,2}(x)\omega_1(x)dxdt\nonumber\\
&\qquad+\int_{t_1}^{T_*}\left(\frac{d}{dt}\frac{2}{2\Omega_1|x_1-x_3|^2}\right)\int f_{2,3}(x)\omega_1(x)dxdt+C\delta T_*^{\epsilon/2-1/3}+I_{2,1}(t_1)\label{eq:ddtI21integrsimp}
\end{align}
We now use use $v_2=O(t^{-1/2})$ on $\supp\omega_1$ and \eqref{eq:dx1dt} (as well as all the analogous statements where we permute the indices $1,2,3$) to get
\begin{equation}\label{cleanupafterIBP}
\frac{d}{dt}\frac{1}{2\Omega_1|x_1-x_2|^2}=O\left(t^{-3/2}\frac{d}{dt}|x_1-x_2|\right)\le O\left(t^{-2}\right)
\end{equation}
as well as all the analogous statements. We plug this along with \eqref{eq:f2bound} into \eqref{eq:ddtI21integrsimp} to get
\begin{align}
I_{2,1}(T_*)&\le CT_*^{-1/2+2\epsilon}(I_{2,1}(t_1)+I_{2,1}(T_*))+C\int_{t_1}^{T_*}\epsilon^2t^{2\epsilon-3/2} I_{k,1}(t)dt+C\delta T_*^{\epsilon/2-1/3}+I_{2,1}(t_1)\nonumber\\
&\le I_{2,1}(t_1)+CT_*^{-1/2+5\epsilon/2}+CT_*^{5\epsilon/2-3/2+2/3}+C\delta T_*^{\epsilon/2-1/3}\nonumber\\
&\le  I_{2,1}(t_1)+C\delta T_*^{\epsilon/2-1/3}\label{I21increment}
\end{align}
where we used $T\gg 1/\delta$. From \eqref{I21increment} and $I_{2,1}(t_1)<t_1^{\epsilon/2}$, we get $I_{2,1}(T_*)<T_*^{\epsilon/2}$.

We now take $k>2$ even, which we treat similarly to the $k=2$ case. We integrate \eqref{eq:ddtIk1compiled} in time from $t_1$ to $T_*$ and use the bootstrap assumptions on $I_{2,1}$ and $I_{k,1}$ to get
\begin{align}
I_{k,1}(T_*)&\le\int_{t_1}^{T_*}Ct^{\epsilon/2}I_{k,1}^{\frac{k-4}{k}}dt-\int_{t_1}^{T_*}\frac{k}{2\Omega_1|x_1-x_2|^2}\frac{d}{dt}\int f_{k,2}(x)\omega_1(x)dxdt\nonumber\\
&\qquad+\int_{t_1}^{T_*}-\frac{k}{2\Omega_1|x_1-x_3|^2}\frac{d}{dt}\int f_{k,3}(x)\omega_1(x)dx\nonumber\\
&\qquad+C\int_{t_1}^{T_*}t^{\frac{k(1+\epsilon)}{4}+\epsilon-\frac{5}{4}}+t^{\frac{(k-1)(1+\epsilon)}{4}+\frac{\epsilon-3}{2}}+\delta t^{\frac{k(1+\epsilon)}{4}-1}+\frac{R^k\delta^{-k}}{t}dt+I_{k,1}(t_1).\nonumber
\end{align}
We integrate by parts in $t$ and use \eqref{cleanupafterIBP} and \eqref{eq:f2bound}, as well as using $T\gg R,1/\delta$ to get
\begin{align}
I_{k,1}(T_*)&\le C\int_{t_1}^{T_*}t^{\epsilon/2}I_{k,1}^{\frac{k-4}{k}}dt+CT_*^{-1/2+2\epsilon}(I_{k,1}(t_1)+I_{k,1}(T_*))+C\int_{t_1}^{T_*}t^{2\epsilon-3/2}I_{k,1}(t)dt\nonumber\\
&\qquad+\int_{t_1}^{T_*}C\delta t^{\frac{k(1+\epsilon)}{4}-1}dt+I_{k,1}(t_1)\nonumber
\end{align}
We now use that $T\gg 1/\delta$, as well as the bootstrap assumption~\ref{bootstrapIki}, to get
\[
I_{k,1}(T_*)\le t_1^{k(1+\epsilon)/4}+CT_*^{-1/2+2\epsilon}T_*^{k(1+\epsilon)/4}+\int_{t_1}^{T_*}C\delta t^{\frac{k(1+\epsilon)}{4}-1}dt
\]
so, as long as $\delta$ is sufficiently small, we get $I_{k,1}(T_*)< T_*^{k(1+\epsilon)/4}$. The same then applies to $I_{k,2}$ and $I_{k,3}$, so bootstrap assumption~\ref{bootstrapIki} is maintained.

We now need to recover bootstrap assumption~\ref{bootstrapsupport}. For this, we let $t_1=T_*-T_*^{2/3}>t_0$ and take some point $p(t_1)$ that is in the support of $\omega_1$, so we have $|p(t_1)-x_1(t_1)|<\epsilon t_1^{1/4+\epsilon}$. We then have $p(t)$ solve $p'(t)=u(t,p(t))$. We want to show that $|p(T_*)-x_1(T_*)|<\epsilon T_*^{1/4+\epsilon}$, which would show that bootstrap assumption~\ref{bootstrapsupport} is maintained since the support of $\omega_1$ is transported by $u$. Suppose this is false, that is $|p(T_*)-x_1(T_*)|\ge \epsilon T_*^{1/4+\epsilon}$. Then let
\[
t_2=\sup\left\{s\in [t_1,T_*]\mid s=t_1\text{ or }|p(s)-x_1(s)|<\frac{\epsilon}{2} T_*^{1/4+\epsilon}\right\}.
\]
 Then we will work on the interval $[t_2,T_*]$, where we are guaranteed that $\frac{\epsilon}{2} T_*^{1/4+\epsilon}\le |p-x_1|\le \epsilon T_*^{1/4+\epsilon}$. We calculate
\begin{align}
\frac{d}{dt}(p-x_1)&=v_2(p)-v_2(x_1)+v_3(p)-v_3(x_1)+\frac{(p-x_1)^\perp}{|p-x_1|^2}\Omega_1+\nonumber\\
&\qquad\qquad+\int \left(\frac{(p-x)^\perp}{|p-x|^2}-\frac{(p-x_1)^\perp}{|p-x_1|^2}\right)\omega_1(x)dx.\label{eq:pprimeexact}
\end{align}
We now use \eqref{eq:v2} to get that
\begin{align}
v_2(p)-v_2(x_1)&=-(p-x_1)\cdot(x_1-x_2)\frac{(x_1-x_2)^\perp}{|x_1-x_2|^4}\nonumber\\
&\qquad-(p-x_1)\cdot(x_1-x_2)^\perp\frac{x_1-x_2}{|x_1-x_2|^4}+O(\epsilon^2 t^{2\epsilon-1})\label{v2errorfine}
\end{align}
We will also need a coarser form of this estimate, namely
\begin{equation}\label{v2errorcoarse}
v_2(p)-v_2(x_1)=O( t^{\epsilon-3/4}).
\end{equation}
We now note that there is some time-dependent matrix $A$ with $|A|\le C|p-x_1|^{-2}$ such that when $|x-x_1|<\frac{|p-x_1|}{2}$, we have that
\begin{equation*}
\frac{(p-x)^\perp}{|p-x|^2}-\frac{(p-x_1)^\perp}{|p-x_1|^2}=A(x-x_1)+O\left(\frac{|x-x_1|^2}{|p-x_1|^3}\right).
\end{equation*}
Let
\begin{align*}
S_7&=\left\{|x-x_1|<\frac{\epsilon}{4} T_*^{1/4+\epsilon}\right\}\\
S_8&=\left\{\epsilon T_*^{1/4+\epsilon}\ge |x-x_1|\ge\frac{\epsilon}{4} T_*^{1/4+\epsilon}\right\}
\end{align*}
Then, since $x_1$ is the center of mass, we have
\begin{align}
&\left|\int \left(\frac{(p-x)^\perp}{|p-x|^2}-\frac{(p-x_1)^\perp}{|p-x_1|^2}\right)\omega_1(x)dx\right|\nonumber\\
&\qquad=\left|\int \frac{(p-x)^\perp}{|p-x|^2}-\frac{(p-x_1)^\perp}{|p-x_1|^2}-A(x-x_1)\omega_1(x)dx\right|\nonumber\\
&\qquad=\int_{S_7} \left|\frac{(p-x)^\perp}{|p-x|^2}-\frac{(p-x_1)^\perp}{|p-x_1|^2}-A(x-x_1)\omega_1(x)\right|dx\nonumber\\
&\qquad\qquad+\int_{S_8} \left|\frac{(p-x)^\perp}{|p-x|^2}-\frac{(p-x_1)^\perp}{|p-x_1|^2}-A(x-x_1)\omega_1(x)\right|dx\nonumber\\
&\qquad\le CI_{2,1}\left(\frac{\epsilon}{2} T_*^{1/4+\epsilon}\right)^{-3}+\int_{O} \frac{1}{|p-x|}+\left(\frac{\epsilon}{2} T_*^{1/4+\epsilon}\right)^{-1}+C|p-x_1|^{-2}|x-x_1|\omega_1(x)dx\nonumber\\
&\qquad\le C T_*^{-3/4-5\epsilon/2}+C\int_{S_8}\frac{1}{|p-x|}\omega_1(x)dx+C T_*^{-1/4-\epsilon}\int_{S_8}\omega_1(x)dx\nonumber\\
&\qquad\qquad+C\left(\frac{\epsilon}{2}T_*^{1/4+\epsilon}\right)^{-2}I_{2,1}^{\half}\left(\int_{S_8}\omega_1(x)dx\right)^{\half}.\label{eq:offcenterpush}
\end{align}
Now, from $I_{k,1}(t)<T_*^{k(1+\epsilon)/4}$, we get that 
\begin{equation}\label{eq:farflungtotalbound}
\int_{S_8}\omega_1(x)dx\le C\frac{T_*^{k(1+\epsilon)/4}}{T_*^{k(1/4+\epsilon)}}=CT_*^{-3k\epsilon/4}.
\end{equation}
From this, we get by Holder's inequality
\begin{equation}\label{eq:farflungimpact}
\int_{S_8}\frac{1}{|p-x|}\omega_1(x)dx\le C\left\|\frac{1}{|p-\cdot|}\right\|_{L^{q_1}(S_8)}\left\|\omega_1\right\|_{L^{2+\epsilon}(S_8)}^{\sigma}\left\|\omega_1\right\|_{L^{1}(S_8)}^{1-\sigma}\le CT_*^{\epsilon-3k\epsilon(1-\sigma)/4}
\end{equation}
for some $\sigma\in (0,1),q_1<2$ that depend only on $\epsilon$. We now choose $k$ sufficiently large that $3k\epsilon(1-\sigma)/4-\epsilon>2$. Then plugging  \eqref{eq:farflungtotalbound} and \eqref{eq:farflungimpact} into \eqref{eq:offcenterpush}, we get that
\begin{equation}\label{offcenterbound}
\int \left(\frac{(p-x)^\perp}{|p-x|^2}-\frac{(p-x_1)^\perp}{|p-x_1|^2}\right)\omega_1(x)dx=O\left(T_*^{-3/4-5\epsilon/2}\right).
\end{equation}
We now plug \eqref{offcenterbound} and \eqref{v2errorcoarse} along with analogous estimate for $v_3$  into \eqref{eq:pprimeexact} to get
\begin{equation}\label{roughpminusxderiv}
\frac{d}{dt}(p-x_1)=\frac{(p-x_1)^\perp}{|p-x_1|^2}\Omega_1+O(T_*^{\epsilon-3/4}).
\end{equation}
Now, for any $\tilde t$ with
\[
[\tilde t, \hat t]:=\left[\tilde t,\tilde t+2\pi|p(\tilde t)-x_1(\tilde t)|^{2}/\Omega_1\right]\subset [t_2,T_*]
\]
we have that on the time interval $[\tilde t, \hat t]$, the total variation of $|p-x_1|$ is at most
\begin{equation}\label{radialvaro}
Var_{[\tilde t, \hat t]}(|p-x_1|)=O(T_*^{\epsilon-3/4}|p(\tilde t)-x_1(\tilde t)|^{2})=O(T_*^{3\epsilon-1/4}).
\end{equation}
We now define $\varphi$ to be the angle of $p-x_1$ and let $\theta_2$ be the angle of $x_1-x_2$. Combining \eqref{roughpminusxderiv} with \eqref{radialvaro} gives that the angular velocity for $t\in [\tilde t,\hat t]$ is
\begin{align*}
\frac{d}{dt}\varphi&=\frac{\Omega_1}{|p-x|^2}+O\left(\frac{T_*^{\epsilon-3/4}}{|p-x|}\right)=\frac{\Omega_1}{|p(\tilde t)-x(\tilde t)|^2}+O\left(\frac{T_*^{3\epsilon-1/4}}{\left(T_*^{1/4+\epsilon}\right)^3}\right)+O(T_*^{-1})\\
&=\frac{\Omega_1}{|p(\tilde t)-x(\tilde t)|^2}+O\left(T_*^{-1}\right)
\end{align*}
so
\begin{equation}\label{angerror}
\varphi(t)=\varphi(\tilde t)+\frac{\Omega_1(t-\tilde t)}{|p(\tilde t)-x(\tilde t)|^2}+O\left(T_*^{2\epsilon-1/2}\right).
\end{equation}
Also, we have $\left|\frac{d}{dt}x_i\right|=O(T_*^{-1/2})$ so
\[
\left|\frac{d}{dt}\theta_2\right|=\left|\frac{d}{dt}\frac{x_1-x_2}{|x_1-x_2|}\right|=O(T_*^{-1})
\]
so
\begin{equation}\label{theta2error}
\theta_2(t)=\theta_2(\tilde t)+O(t_1^{2\epsilon-1/2})
\end{equation}
We now use \eqref{eq:pprimeexact} and \eqref{v2errorfine} along with the analogous estimate for $x_3$ and \eqref{offcenterbound} to compute
\begin{align}
\frac{d}{dt}|p-x_1|^2&=-4(p-x_1)\cdot(x_1-x_2)\frac{(p-x_1)\cdot(x_1-x_2)^\perp}{|x_1-x_2|^4}\nonumber\\
&\qquad-4(p-x_1)\cdot(x_1-x_3)\frac{(p-x_1)\cdot(x_1-x_3)^\perp}{|x_1-x_3|^4}\nonumber\\
&\qquad+O(T_*^{(2\epsilon-1)+(1/4+\epsilon)})+O\left(T_*^{(-3/4-5\epsilon/2)+(1/4+\epsilon)}\right)\nonumber\\
&=\frac{-2\sin(2\theta_2+2\varphi)|p-x_1|^2}{|x_1-x_2|^2}-\frac{2\sin(2\theta_3+2\varphi)|p-x_1|^2}{|x_1-x_3|^2}+O\left(T_*^{-1/2-3\epsilon/2}\right).\label{eq:squareddistderiv}
\end{align}
By combining \eqref{theta2error}, \eqref{angerror}, \eqref{radialvaro}, and \eqref{cleanupafterIBP}, we get that for $t\in[\tilde t,\hat t]$, we have
\begin{align}
&\frac{-2\sin(2\theta_2+2\varphi)|p-x_1|^2}{|x_1-x_2|^2}=\nonumber\\
&\qquad\qquad=-2\sin\left(2\theta_2(\tilde t)+2\varphi(\tilde t)+\frac{2\Omega_1(t-\tilde t)}{|p(\tilde t)-x(\tilde t)|^2}\right)\frac{|p(\tilde t)-x_1(\tilde t)|^2}{|(x_1(\tilde t)-x_2(\tilde t))|^2}+O\left(T_*^{4\epsilon-1}\right).\label{sin2error}
\end{align}
Then substituting \eqref{sin2error} into \eqref{eq:squareddistderiv} (along with the analogous estimate for $x_3$ and integrating from $\tilde t$ to $\hat t$, we note that the principal term of \eqref{sin2error} cancels and we are left with
\[
|p(\hat t)-x_1(\hat t)|^2-|p(\tilde t)-x_1(\tilde t)|^2=O((\hat t-\tilde t)T_*^{-3\epsilon/2-1/2}).
\]
We now take a new interval starting at $\hat t$. Tiling most of $[t_2,T_*]$ with such intervals, we get that
\begin{equation}\label{pminusxsquared}
|p(T_*)-x_1(T_*)|^2\le |p(t_2)-x_1(t_2)|^2+O\left((T_*-t_2)T_*^{-3\epsilon/2-1/2}\right)+\int_{t_3}^{T_*}\frac{d}{dt}|p-x_1|^2dt
\end{equation}
where $t_3\in [t_2,T_*]$ satisfies $T_*-t_3=O(t^{1/2+2\epsilon})$. Using \eqref{roughpminusxderiv} to bound the last term of \eqref{pminusxsquared}, we then get
\[
|p(T_*)-x_1(T_*)|^2\le |p(t_1)-x_1(t_1)|^2+O(T_*^{1/6-3\epsilon/2})+O(T_*^{4\epsilon})\le(\epsilon t_1^{1/4+\epsilon})^2+O(T_*^{1/6-3\epsilon/2})<\left(\epsilon T_*^{1/4+\epsilon}\right)^2,
\]
which verifies bootstrap assumption~\ref{bootstrapsupport}. Thus we have shown (modulo the proof of \eqref{eq:subsecresult} in section~\ref{derivest}) that we cannot have $t_0+t_0^{9/10}\le T_*<\infty$.

\subsection{Moment renormalization}\label{derivest}
In this section, we prove estimate~\eqref{eq:subsecresult}. This estimate would follow from a short computation directly if all of the mass of the vortex patch were located precisely at $x_i$ and all of the $k$th moment came from parts of the vortex patch at a large distance from $x_i$. This cannot hold precisely, but we will obtain an approximation to this by proving that each $\omega_i$ concentrates, as shown in Figure 2 (see \eqref{concentration} below for the precise statement).
\begin{figure}[h!]
\centering
\includegraphics[scale=0.5]{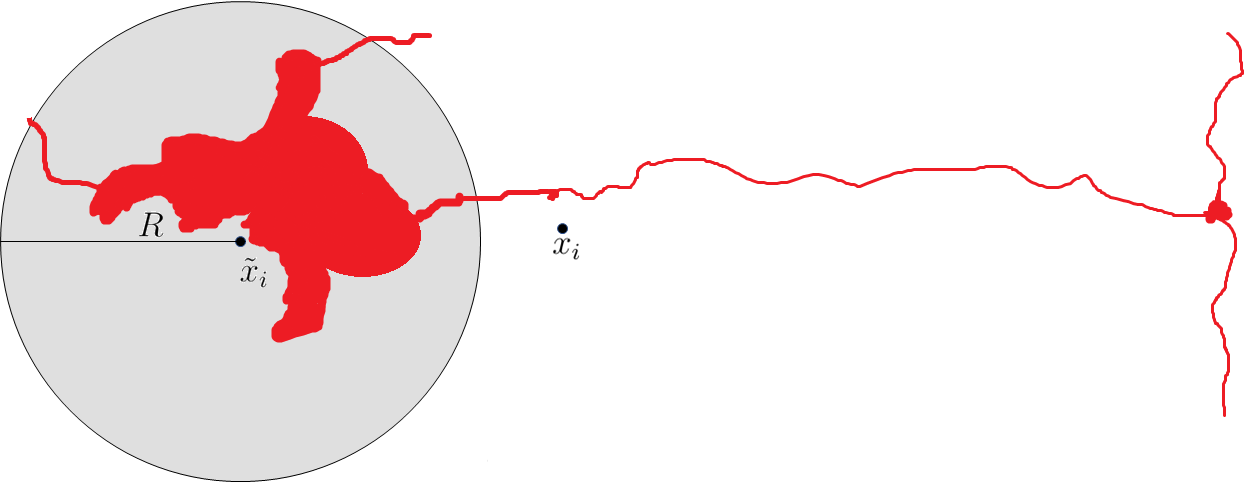}
\caption{The red mass is the vortex patch given by $\omega_i$. Note that most of its mass is in a ball centered at some $\tilde x_i$ that is distinct from the center of mass $x_i$.}
\end{figure}

To prove the concentration result, we note that for any solution of 2D Euler with compactly supported $L^{2+\epsilon}$ vorticity, the following is a conserved quantity:
\[
L=\iint \log|x-y| \omega(x)\omega(y)dxdy
\]
This quantity corresponds to physical energy of the fluid, and one can directly check that it is conserved with a simple computation. For our solution, we have that for any two points $x\in \supp\omega_i,y\in\supp\omega_j$ with $i\ne j$, we have $\log|x-y|=(\log t)/2+O(1)$. Thus
\bal
L&= \sum_{i\ne j}\Omega_i\Omega_j(\log t)/2+O(1)+\sum_{i=1}^3\iint\omega_i(x)\omega_i(y)\log|x-y|  dxdy\\
&=\sum_{i=1}^3\iint\omega_i(x)\omega_i(y)\log|x-y|  dxdy+O(1)
\eal
where we used \eqref{eq:cond1}. Now note that $L$ is conserved, and using a combination of Holder's and Young's inequality, we get
\[
\left|\iint\omega_i(x)\omega_i(y)(\log|x-y|)_-  dxdy\right|\le \|\omega_i\|_{L^{2+\epsilon}}^2\|\log|\cdot|\|_{L^{1/(2-2/(2+\epsilon))}}=O(1).
\]
Thus there is some constant $C$ with
\[
\iint\omega_i(x)\omega_i(y)(\log|x-y|)_{+}  dxdy\le C.
\]
Thus there is some $R=R(\delta)>0$ such that 
\[
\iint\omega_i(x)\omega_i(y)\mathbbm{1}_{|x-y|>R}  dxdy<\Omega_i\delta^4,
\]
from which it follows that for some $\tilde x_i\in\supp \omega_i$, we have that the vorticity mass of $\omega_i$ concentrates around $\tilde x_i$, meaning that
\begin{equation}\label{concentration}
\int\omega_i(x)\mathbbm{1}_{|x-\tilde x_i|>R} dx<\delta^4.
\end{equation}
We now recall that we defined $f_{k,2}(x)=-\cos(2\theta(x))|x-x_1|^{k+2}$. In the calculation below, we will use the following facts:
\bal
|f|&\le |x-x_1|^{k+2}\\
|Df_{k,2}|&\le C|x-x_1|^{k+1}\\
|D^2f_{k,2}|&\le C|x-x_1|^{k}
\eal
We calculate
\begin{align}
\frac{d}{dt} &\int f_{k,2}(x)\omega_1(x)dx=\nonumber\\
&=\iint Df_{k,2}(x)\cdot\frac{(x-y)^\perp}{|x-y|^2}\omega_1(x)\omega_1(y)dxdy+\int Df_{k,2}(x)\cdot(v_2(x)+v_3(x))\omega_1(x)dx.\label{eq:fderiv}
\end{align}
In the calculation below, we will use
\begin{align*}
S_4&=\{|x-x_1|\le \delta|y-x_1|\}\\
S_5&=\{\delta|y-x_1|<|x-x_1|< |y-x_1|/\delta\}\\
S_6&=\{|x-x_1|\ge|y-x_1|/\delta\}.
\end{align*}
Using $v_2,v_3=O(t^{-1/2})$, and symmetrizing in $x$ and $y$ for the integral over $S_5$ we get 
\begin{align}
\frac{d}{dt}& \int f_{k,2}(x)\omega_1(x)dx\nonumber\\
&=\iint_{S_4\cup S_5\cup S_6} Df_{k,2}(x)\cdot\frac{(x-y)^\perp}{|x-y|^2}\omega_1(x)\omega_1(y)dxdy+ O\left(t^{-1/2}\right)\int|x-x_1|^{k+1}\omega_1(x)dx\nonumber\\
&=\iint_{S_4}Df_{k,2}(x)\cdot\frac{(x-y)^\perp}{|x-y|^2}\omega_1(x)\omega_1(y)dxdy\nonumber\\
&\qquad+\half\iint_{S_5}(Df_{k,2}(x)-Df_{k,2}(y))\cdot\frac{(x-y)^\perp}{|x-y|^2}\omega_1(x)\omega_1(y)dxdy\nonumber\\
&\qquad+\iint_{S_6}Df_{k,2}(x)\cdot\frac{(x-y)^\perp}{|x-y|^2}-Df_{k,2}(x)\cdot\frac{(x-x_1)^\perp}{|x-x_1|^2}\omega_1(x)\omega_1(y)dxdy\nonumber\\
&\qquad +\iint_{S_4\cup S_5}- Df_{k,2}(x)\cdot\frac{(x-x_1)^\perp}{|x-x_1|^2}\omega_1(x)\omega_1(y)dxdy\nonumber\\
&\qquad +\iint_{\R^3\times\R^3}Df_{k,2}(x)\cdot\frac{(x-x_1)^\perp}{|x-x_1|^2}\omega_1(x)\omega_1(y)dxdy\nonumber\\
&\qquad+ O\left(t^{-1/2}\sup_{x\in \supp\omega_1}|x-x_1|\int|x-x_1|^{k}\omega_1(x)dx\right).\label{eq:derivativemess}
\end{align}
We now  address each of these terms. First, we have (one needs to be careful about the sign of $\theta$ here)
\begin{equation}\label{eq:principalterm}
\iint_{\R^3\times\R^3} Df_{k,2}(x)\cdot\frac{(x-x_1)^\perp}{|x-x_1|^2}\omega_1(x)\omega_1(y) dxdy=\Omega_1\int 2\sin(2\theta)|x-x_1|^k\omega_1(x) dx.
\end{equation}
This is the principal term. The rest we need to bound with something small. The first term we bound by
\begin{align}
&\left|\iint_{S_4}Df_{k,2}(x)\cdot\frac{(x-y)^\perp}{|x-y|^2}\omega_1(x)\omega_1(y)dxdy\right|\nonumber\\
&\qquad\le C\iint_{S_4}\frac{|x-x_1|^{k+1}}{|y-x_1|}\omega_1(x)\omega_1(y)dxdy\nonumber\\
&\qquad\le C \delta\iint_{\R^3\times\R^3} |x-x_1|^{k}\omega_1(x)\omega_1(y)dxdy\nonumber\\
&\qquad\le \delta C I_{k,1}.\label{eq:firstjunk}
\end{align}
We bound the third term in \eqref{eq:derivativemess} by
\begin{align}
&\left|\iint_{S_6}Df_{k,2}(x)\cdot\frac{(x-y)^\perp}{|x-y|^2}-Df_{k,2}(x)\cdot\frac{(x-x_1)^\perp}{|x-x_1|^2}\omega_1(x)\omega_1(y)dxd\right|\nonumber\\
&\qquad\le C\iint_{S_5}|x-x_1|^{k+1}\frac{|y-x_1|}{|x-x_1|^2}\omega_1(x)\omega_1(y)dxdy\nonumber\\
&\qquad\le \delta C\iint_{S_4\cup S_5}|x-x_1|^{k}\omega_1(x)\omega_1(y)dxdy\nonumber\\
&\qquad\le \delta CI_{k,1}.\label{eq:thirdjunk}
\end{align}
 We bound the second term in \eqref{eq:derivativemess} by using the mean value theorem on $Df$ to get
\begin{align}
&\left|\iint_{S_5}(Df_{k,2}(x)-Df_{k,2}(y))\cdot\frac{(x-y)^\perp}{|x-y|^2}\omega_1(x)\omega_1(y)dxdy\right|\nonumber\\
&\qquad\le C\iint_{S_5}(|x-x_1|^{k}+|y-x_1|^{k})\frac{|x-y|^2}{|x-y|^2}\omega_1(x)\omega_1(y)dxdy\nonumber\\
&\qquad\le C\iint_{S_5}|x-x_1|^{k}\omega_1(x)\omega_1(y)dxdy\nonumber\\
&\qquad\le C\iint_{S_4\cup S_5}|x-x_1|^{k}\omega_1(x)\omega_1(y)dxdy.\label{eq:secondjunk}
\end{align}
We bound the fourth term in \eqref{eq:derivativemess} by
\begin{align}
&\left|\iint_{S_4\cup S_5}- Df_{k,2}(x)\cdot\frac{(x-x_1)^\perp}{|x-x_1|^2}\omega_1(x)\omega_1(y)dxdy\right|\nonumber\\
&\qquad\le C\iint_{S_4\cup S_5}|x-x_1|^{k}\omega_1(x)\omega_1(y)dxdy.\label{eq:fourthjunk}
\end{align}
To bound this expression, we use the concentration result \eqref{concentration} and split into two cases. First, if $|\tilde x_1-x_1|<2R$, then
\begin{align}
\iint_{S_4\cup S_5}&|x-x_1|^k\omega_1(x)\omega_1(y)dxdy\nonumber\\
&\le \iint_{\{|y-x_1|\ge 3R\}}|x-x_1|^k\omega_1(x)\omega_1(y)dxdy+\iint_{\{|x-x_1|\le 3R/\delta\}}|x-x_1|^k\omega_1(x)\omega_1(y)dxdy\nonumber\\
&\le \delta^4 \Omega_1I_{k,1}+\Omega_1^2(3R/\delta)^k.\label{eq:caseonebound}
\end{align}
The second case is $|\tilde x_1-x_1|\ge 2R$. In this case, we have (since $\delta$ is sufficiently small)
\[
\left|\int_{B(\tilde x_1,R)}(x-x_1)\omega_1(x)dx\right|\ge \frac{|\tilde x_1-x_1|\Omega_1}{4}.
\]
Also,
\[
\int (x-x_1)\omega_1(x) dx=0
\]
and by  \eqref{concentration}, we have
\[
\int_{B(0,|\tilde x_1-x_1|\Omega_1/(8\delta^4))\wo B(\tilde x_1,R)}|x-x_1|\omega_1(x)dx\le \frac{|\tilde x_1-x_1|\Omega_1}{8}.
\]
Combining the last three inequalities, we get
\[
\int_{\{|x-x_1|>|\tilde x_1-x_1|\Omega_1/(8\delta^4)\}}|x-x_1|\omega_1(x)dx\ge \frac{|\tilde x_1-x_1|\Omega_1}{8}
\]
from which it follows that
\[
I_{k,1}\ge \int_{\{|x-x_1|>|\tilde x_1-x_1|\Omega_1/(8\delta^4)\}}|x-x_1|^k\omega_1(x)dx\ge  \frac{|\tilde x_1-x_1|^k\Omega_1^k}{8^k\delta^{4(k-1)}}.
\]
Thus
\[
\int_{\{|x-x_1|\le 2|\tilde x_1-x_1|/\delta\}}|x-x_1|^k\omega_1(x)dx\le  2^k\Omega_1|\tilde x_1-x_1|^k/\delta^k<\delta I_{k,1}.
\]
as long as $\delta$ is sufficiently small. Thus
\begin{align}
\iint_{S_4\cup S_5}&|x-x_1|^{k}\omega_1(x)\omega_1(y)dxdy\nonumber\\
&\le \iint_{\{|x-x_1|<2|\tilde x_1-x_1|/\delta\}}|x-x_1|^k\omega_1(x)\omega_1(y)dxdy\nonumber\\
&+\iint_{\{|y-x_1|\ge 2|\tilde x_1-x_1|\}}|x-x_1|^k\omega_1(x)\omega_1(y)dxdy\nonumber\\
&\le \delta I_{k,1}+\delta^4 I_{k,1}=2\delta I_{k,1}.\label{eq:casetwobound}
\end{align}
Thus, combining \eqref{eq:caseonebound} and \eqref{eq:casetwobound}, we get that in either case, we have
\begin{equation}\label{eq:sparsekthmoment}
\iint_{S_4\cup S_5}|x-x_1|^k\omega_1(x)\omega_1(y)dxdy\le\delta  CI_{k,1}+C\delta^{-k}R^k.
\end{equation}
Combining \eqref{eq:derivativemess}, with \eqref{eq:principalterm}, and getting error bounds from \eqref{eq:firstjunk},\eqref{eq:thirdjunk},\eqref{eq:secondjunk},\eqref{eq:fourthjunk}, and \eqref{eq:sparsekthmoment}, we get 
\begin{equation*}
\frac{d}{dt}\int f_{k,2}(x)\omega_1(x)dx=\Omega_1\int 2\sin(2\theta)|x-x_1|^k\omega_1(x) dx+O(\delta I_{k,1}+R^k\delta^{-k}).
\end{equation*}
which completes the proof of \eqref{eq:subsecresult}.
\section{Short time behavior}
In this section, we assume that $T_*<t_0+t_0^{9/10}$. We will use rougher versions of estimates from the section~\ref{sec:longterm}. In particular, we use the boundedness of $\sin$ to turn \eqref{eq:dotv2} into
\begin{equation*}
\int |x-x_1|^{k-2}(x-x_1)\cdot v_2(x)dx=O\left(t^{-1}I_{k,1}\right)+O\left(t^{\epsilon-5/4}I_{k,1}\right)+O\left(t^{\frac{\epsilon-3}{2}}I_{k,1}^{\frac{k-1}{k}}\right).
\end{equation*}
We then plug this, the analogous bound for $v_3$, and \eqref{usualmomentderiv} into \eqref{eq:ddtIk1} to get that whenever $I_{k,1}\ge 1$, we have
\begin{equation}\label{earlymomentderiv}
\frac{d}{dt}I_{k,1}\le C\frac{I_{k,1}}{t}+ C(k-2)t^{\epsilon/2}I_{k,1}^{\frac{k-4}{k}}.
\end{equation}
Plugging in $k=2$, and using the fact that $I_{2,1}(t_0)=O(1)$, we get that for all $t\in [t_0,T_*]$, we have
\[
I_{2,1}(t)=O(1)\exp\left(\int_{t_0}^{t_0+t_0^{9/10}}\frac{C}{t}dt\right)=O(1)<t^{\epsilon/2}
\]
so bootstrap assumption~\ref{bootstrapI2i} is maintained.
Now we apply \eqref{earlymomentderiv} for more general $k$. Using bootstrap assumption~\ref{bootstrapIki}, we get that
\begin{align*}
I_{k,1}(T_*)&\le C+C\int_{t_0}^{T_*}\frac{I_{k,1}}{t}+t^{\epsilon/2}I_{k,1}^{\frac{k-4}{k}}dt\\
                    &\le C+Ct_0^{9/10}\left(t_0^{k(1+\epsilon)/4-1}+t_0^{\epsilon/2}t_0^{k(1+\epsilon)/4-(1+\epsilon)}\right)\\
                    &<T_*^{k(1+\epsilon)/4}
\end{align*}
so bootstrap assumption~\ref{bootstrapIki} is maintained. Now, to verify bootstrap assumption~\ref{bootstrapsupport}, we do things similarly to section~\ref{sec:longterm}. We suppose that there is some $p(t)\in\supp\omega_1$ transported by $u$ such that at time $T_*$, we have $|p(T_*)-x_1(T_*)|\ge \epsilon T_*^{1/4+\epsilon}$ and we define
\[
t_2=\sup\left\{s\in [t_0,T_*]\,\middle\vert\, s=t_0\text{ or }|p(s)-x_1(s)|<\frac{\epsilon}{2} T_*^{1/4+\epsilon}\right\}.
\]
Since $|p(t_0)-x_1(t_0)|\le \rho$, we in fact have that $t_2>t_0$ and that
\[
 |p(t_2)-x(t_2)|=\frac{\epsilon}{2} T_*^{1/4+\epsilon}.
\]
Then on the interval $[t_2,T_*]$, we have that \eqref{roughpminusxderiv} holds, so
\[
|p(T_*)-x(T_*)|\le |p(t_2)-x(t_2)|+O\left(T_*^{9/10}T_*^{\epsilon-3/4}\right)< \epsilon T_*^{1/4+\epsilon},
\]
which gives a contradiction, so bootstrap assumption~\ref{bootstrapsupport} is maintained. Thus we have shown that we cannot have $T_*<t_0+t_0^{9/10}$, so $T_*=\infty$, completing the proof of Theorem~\ref{thm:threepatches}.

\appendix
\section{Proof of Lemma}
Here we prove Lemma~\ref{configlemma}. For the first part of the statement, we note that by \eqref{eq:cond1},
\[
\left(\sum_{i=1}^3 \Omega_i\right)^2=\sum_{i=1}^3 \Omega_i^2+2\sum_{i<j}\Omega_i\Omega_j>0
\]
so $\sum\Omega_i\ne 0$. For the second part of the lemma statement we first want to show that $\nabla E$ and $\nabla I$ are not parallel. We define
\bal
\tilde I&=\sum_{i=1}^3\sum_{j=1}^3 \Omega_i\Omega_j|x_i-x_j|^2=-2\left|\sum_{i=1}^3\Omega_i x_i\right|^2+2\sum_{i=1}^3\Omega_i\left(\sum_{j=1}^3 \Omega_j\right)|x_i|^2\\
&=-2X^2+2\left(\sum_{j=1}^3 \Omega_j\right)I=2\left(\sum_{j=1}^3 \Omega_j\right)I
\eal
so, on the space $V$, we have that $\tilde I$ is proportional to $I$. Now, we can think of any small perturbation as a small change in $|x_1-x_2|,|x_1-x_3|,|x_2-x_3|$, plus a rotation (this works because the points are not collinear, so there is some room in the triangle inequality). Then
\[
\frac{\partial E}{\partial |x_i-x_j|}=\frac{\Omega_i\Omega_j}{|x_i-x_j|},\frac{\partial \tilde I}{\partial |x_i-x_j|}=2 \Omega_i\Omega_j|x_i-x_j|.
\]
Then the only way $\nabla E$ and $\nabla\tilde I$ are parallel is if all $|x_i-x_j|$ are equal, that is the points form the vertices of an equilateral triangle. The condition of the lemma excludes this case specifically. Now, since the gradients are not parallel, we can locally find a surface $S$ in $A$ through $(x_1,x_2,x_3)$ on which $E$ and $\tilde I$ give coordinates. Since $\tilde I$ is proportional to $I$, we have that $E$ and $I$ give coordinates. Rotation clearly does not change the values of $E$ or $I$. Also, at the point $(x_1,x_2,x_3)$, we have by conditions \eqref{eq:cond1} and \eqref{eq:cond2} that scaling also does not change the values of $I$ and $E$. Thus, rotation and scaling give two vectors fields that are linearly independent and whose span does not intersect the tangent space of $S$. Thus, in some small set $U$, we can take the coordinates given in the lemma statement.
\bibliographystyle{abbrv}
\bibliography{../zbarskybib}

\begin{thebibliography}{10}

\bibitem{Aref07}
H.~Aref.
\newblock Point vortex dynamics: a classical mathematics playground.
\newblock {\em J. Math. Phys.}, 48(6):065401, 23, 2007.

\bibitem{ConstantinBertozzi93}
A.~L. Bertozzi and P.~Constantin.
\newblock Global regularity for vortex patches.
\newblock {\em Comm. Math. Phys.}, 152(1):19--28, 1993.

\bibitem{Burbea82}
J.~Burbea.
\newblock Motions of vortex patches.
\newblock {\em Lett. Math. Phys.}, 6(1):1--16, 1982.

\bibitem{Chemin91}
J.-Y. Chemin.
\newblock Persistance des structures g\'{e}om\'{e}triques li\'{e}es aux poches
  de tourbillon.
\newblock In {\em S\'{e}minaire sur les \'{E}quations aux {D}\'{e}riv\'{e}es
  {P}artielles, 1990--1991}, pages Exp. No. XIII, 11. \'{E}cole Polytech.,
  Palaiseau, 1991.

\bibitem{ChoiDenisov19}
K.~Choi and S.~Denisov.
\newblock On the growth of the support of positive vorticity for 2{D} {E}uler
  equation in an infinite cylinder.
\newblock {\em Comm. Math. Phys.}, 367(3):1077--1093, 2019.

\bibitem{DeemZabusky78}
G.~S. {Deem} and N.~J. {Zabusky}.
\newblock {Vortex Waves: Stationary ``V States,'' Interactions, Recurrence, and
  Breaking.}
\newblock {\em Phys. Rev. Lett.}, 41(7):518, Aug 1978.

\bibitem{DurrPulvurenti82}
D.~D\"{u}rr and M.~Pulvirenti.
\newblock On the vortex flow in bounded domains.
\newblock {\em Comm. Math. Phys.}, 85(2):265--273, 1982.

\bibitem{ElgindiJeong19}
T.~Elgindi and I.-J. Jeong.
\newblock On singular vortex patches, ii: Long-time dynamics.
\newblock 2019.
\newblock preprint, \url{https://arxiv.org/abs/1909.13555}.

\bibitem{GomezSerranoetal19}
J.~Gomez-Serrano, J.~Park, J.~Shi, and Y.~Yao.
\newblock Symmetry in stationary and uniformly rotating solutions of active
  scalar equations.
\newblock 2019.
\newblock preprint, \url{https://arxiv.org/abs/1908.01722}.

\bibitem{Iftimie99}
D.~Iftimie.
\newblock \'{E}volution de tourbillon \`a support compact.
\newblock In {\em Journ\'{e}es ``\'{E}quations aux {D}\'{e}riv\'{e}es
  {P}artielles'' ({S}aint-{J}ean-de-{M}onts, 1999)}, pages Exp. No. IV, 8.
  Univ. Nantes, Nantes, 1999.

\bibitem{Iftimie07}
D.~Iftimie.
\newblock Large time behavior in perfect incompressible flows.
\newblock In {\em Partial differential equations and applications}, volume~15
  of {\em S\'{e}min. Congr.}, pages 119--179. Soc. Math. France, Paris, 2007.

\bibitem{SecondIftimieetal03}
D.~Iftimie, M.~C. Lopes~Filho, and H.~J. Nussenzveig~Lopes.
\newblock Large time behavior for vortex evolution in the half-plane.
\newblock {\em Comm. Math. Phys.}, 237(3):441--469, 2003.

\bibitem{Iftimieetal03}
D.~Iftimie, M.~C. Lopes~Filho, and H.~J. Nussenzveig~Lopes.
\newblock On the large-time behavior of two-dimensional vortex dynamics.
\newblock {\em Phys. D}, 179(3-4):153--160, 2003.

\bibitem{Iftimieetal07}
D.~Iftimie, M.~C. Lopes~Filho, and H.~J. Nussenzveig~Lopes.
\newblock Confinement of vorticity in two dimensional ideal incompressible
  exterior flow.
\newblock {\em Quart. Appl. Math.}, 65(3):499--521, 2007.

\bibitem{IftimieMarchioro18}
D.~Iftimie and C.~Marchioro.
\newblock Self-similar point vortices and confinement of vorticity.
\newblock {\em Comm. Partial Differential Equations}, 43(3):347--363, 2018.

\bibitem{InftimieSiderisGamblin99}
D.~Iftimie, T.~C. Sideris, and P.~Gamblin.
\newblock On the evolution of compactly supported planar vorticity.
\newblock {\em Comm. Partial Differential Equations}, 24(9-10):1709--1730,
  1999.

\bibitem{Kirchoff1874}
G.~Kirchhoff.
\newblock {\em Vorlesungen \:{u}ber mathematische {P}hysik}.
\newblock Teubner, 1874.

\bibitem{Kudela14}
H.~Kudela.
\newblock Collapse of {$n$}-point vortices in self-similar motion.
\newblock {\em Fluid Dyn. Res.}, 46(3):031414, 16, 2014.

\bibitem{LopesLopes98}
M.~C. Lopes~Filho and H.~J. Nussenzveig~Lopes.
\newblock An extension of {M}archioro's bound on the growth of a vortex patch
  to flows with {$L^p$} vorticity.
\newblock {\em SIAM J. Math. Anal.}, 29(3):596--599, 1998.

\bibitem{Marchioro94}
C.~Marchioro.
\newblock Bounds on the growth of the support of a vortex patch.
\newblock {\em Comm. Math. Phys.}, 164(3):507--524, 1994.

\bibitem{Marchioro96}
C.~Marchioro.
\newblock On the growth of the vorticity support for an incompressible
  non-viscous fluid in a two-dimensional exterior domain.
\newblock {\em Math. Methods Appl. Sci.}, 19(1):53--62, 1996.

\bibitem{Marchioro98}
C.~Marchioro.
\newblock On the localization of the vortices.
\newblock {\em Boll. Unione Mat. Ital. Sez. B Artic. Ric. Mat. (8)},
  1(3):571--584, 1998.

\bibitem{MarchioroPulvurenti93}
C.~Marchioro and M.~Pulvirenti.
\newblock Vortices and localization in {E}uler flows.
\newblock {\em Comm. Math. Phys.}, 154(1):49--61, 1993.

\bibitem{MarchioroPulvurenti94}
C.~Marchioro and M.~Pulvirenti.
\newblock {\em Mathematical theory of incompressible nonviscous fluids},
  volume~96 of {\em Applied Mathematical Sciences}.
\newblock Springer-Verlag, New York, 1994.

\bibitem{NovikovSedov79}
E.~Novikov and Y.~Sedov.
\newblock Vortex collapse.
\newblock {\em J. Exp. Theor. Phys.}, (50):297–301, 1979.

\bibitem{Serfati94}
P.~Serfati.
\newblock Une preuve directe d'existence globale des vortex patches {$2$}{D}.
\newblock {\em C. R. Acad. Sci. Paris S\'{e}r. I Math.}, 318(6):515--518, 1994.

\bibitem{Serfati98}
P.~Serfati.
\newblock Borne en temps des caract\'eristiques de l'\'equation d'{E}uler 2d
  \`a tourbillon positif et localisation pour le mod\`ele point-vortex, 1998.
\newblock Manuscript.

\bibitem{OtherSerfati98}
P.~Serfati.
\newblock Tourbillons-presque-mesures spatialement born\'{e}s et \'{e}quation
  d'{E}uler 2{D}, 1998.
\newblock Manuscript.

\bibitem{Yudovich63}
V.~I. Yudovich.
\newblock Non-stationary flows of an ideal incompressible fluid.
\newblock {\em \v{Z}. Vy\v{c}isl. Mat i Mat. Fiz.}, 3:1032--1066, 1963.

\end{thebibliography}
\end{document}